\documentclass[reqno]{amsart}
\usepackage{url}
\usepackage{amssymb}
\usepackage{graphicx}
\usepackage[colorinlistoftodos]{todonotes}
\usepackage{amsmath}
\usepackage{amsthm}

\newtheorem{theorem}{Theorem}[section]
\newtheorem{lemma}[theorem]{Lemma}
\newtheorem{corollary}[theorem]{Corollary}

\theoremstyle{definition}

\theoremstyle{remark}

\numberwithin{equation}{section}
\usepackage{algorithmic}
\usepackage[top=3cm,bottom=2cm,right=2cm,left=2cm]{geometry}

\begin{document}

\title[On the Classification of Modular Congruence Families]{On the Classification of Modular Congruence Families}


\author{}
\address{}
\curraddr{}
\email{}
\thanks{}

\author{Nicolas Allen Smoot}
\address{}
\curraddr{}
\email{}
\thanks{}



\date{}

\dedicatory{}

\begin{abstract}
Congruence families, i.e., $\ell$-adic convergence for well-defined arithmetic subsequences, is a commonplace phenomenon for the coefficients of modular forms.  Such families superficially resemble one another, but they often vary substantially in difficulty.  Moreover, the critical difficulties associated with a given family will generally manifest themselves at the end stages of an attempted proof.  We give a conjectured classification system of congruence families for the coefficients of modular eta quotients by studying the topology of the associated modular curve.
\end{abstract}

\maketitle

\section{Introduction}

\subsection{Background}

The theory of integer partitions has a long and rich history, having been studied in some form since Leibniz, and having been systematically studied since Euler \cite{Euler}.  If we define $p(n)$ as the counting function for the number of unrestricted partitions of $n\in\mathbb{N}_0$, then the sequence for $p(n)$ begins

\begin{align*}
\left(p(n)\right)_{n\ge 0} =\big(&1, 1, 2, 3, 5, 7, 11, 15, 22, 30, 42, 56, 77, 101, 135, 176, 231, 297, 385, 490, 627, 792, 1002, 1255, 1575, 1958,...\big).
\end{align*}  Euler and his contemporaries appear to have regarded the sequence as pseudorandom, and there is no indication that any substantial arithmetic properties for $p(n)$ were found before the twentieth century.  Indeed, the sequence for $p(n)$ has been compared with that of the prime counting function $\pi(x)$ with respect to its apparent lack of predictability.

This perspective was radically altered almost overnight in 1917-18 when Ramanujan discovered some astonishing arithmetic structure hiding in the apparent randomness of the sequence---all the more remarkable given how accessible this structure appears \cite{Ramanujan}.  Simply starting with the fifth term in the sequence above ($p(4)=5$) and counting every fifth term thereafter, we have the subsequence

\begin{align*}
\left(p(5n+4)\right)_{n\ge 0} =\big(&5, 30, 135, 490, 1575,...\big).
\end{align*}  The unexpected divisibility of $p(5n+4)$ by 5---which Ramanujan immediately proved---is a celebrated achievement in the subject.  Moreover, once Ramanujan knew what sort of patterns to look for, he not only found analogous properties for $p(n)$ modulo 7 and 11, but also for $5^2$, $7^2$, and beyond.  He quickly organized his results into a sweeping conjecture that for $\ell\in\{5,7,11\}$, $p(n)$ must be divisible by $\ell^{\alpha}$ for any $\alpha\ge 1$ whenever $24n-1$ is divisible by $\ell^{\alpha}$ \cite{Ramanujan}.

This conjecture led to half a century of work in the subject.  The case for $\ell=5$ was first published in 1938 by Watson \cite{Watson} (although Ramanujan appears to have understood the proof some twenty years prior \cite{BO}).  The case for $\ell=7$ was found to be more interesting than Ramanujan had anticipated, and Watson published the proof of a slight modification of Ramanujan's conjecture in the same year \cite{Watson}.

On the other hand, the case for $\ell=11$ is not easy at all.  An acceptable proof was not provided until 1967 \cite{Atkin}.  The complete theorem stands thus:

\begin{theorem}\label{ramp5711}
\begin{align}
\text{If } n,\alpha\in\mathbb{Z}_{\ge 1} \text{ and } 24n\equiv 1\pmod{\ell^{\alpha}},&\text{ then } p(n)\equiv 0\pmod{\ell^{\beta}},\label{origconjecR5711}
\end{align} with

\begin{align}
\beta := 
\begin{cases}
\alpha & \text{ if } \ell\in\{5,11\},\\
\left\lfloor \frac{\alpha}{2} \right\rfloor + 1 & \text{ if } \ell=7.
\end{cases}\label{origconjecR5711B}
\end{align}
\end{theorem}

On its own, this theorem is remarkable.  Still more striking is the fact that similar behavior has been found to occur not only in $p(n)$, but also in a very large class of integer sequences.

For example, the counting function for partitions into distinct parts, here denoted $p_D(n)$, exhibits the following congruence family \cite{Rodseth}:

\begin{theorem}[R\o dseth]\label{distinctpartsthm}
For $\alpha\ge 1$ and $n\ge 0$, the following applies:
\begin{align}
\text{If } 24n\equiv -1\pmod{5^{2\alpha+1}}&\text{, then } p_{D}(n)\equiv 0\pmod{5^{\alpha}}.
\end{align}
\end{theorem}  This has a very similar appearance to the families for $p(n)$.  Similar examples can be found to apply to counting functions for $k$-colored partitions \cite{Atkin0}, $k$-elongated plane partition diamonds \cite{Smoot}, \cite{SellersB}, \cite{Baner5}, \cite{Baner7}, Frobenius partitions \cite{Paule}, \cite{SellersA}, and many others.

In the last few years I have labored to create a systematic approach to congruence families.  Certainly, many of the initial steps in the proofs are already very standardized.  However, there is substantial variation in the later steps, and no one approach will suffice for all problems.

In trying to unify the various different methods I have discovered a pattern in the general difficulty and accessibility of a given congruence family based on its associated modular curve.

One of the most noteworthy problems with proving congruence families is that the most serious difficulties generally arise in the final steps of the proof.  This means that one can commit a very large degree of work to a proof which looks extremely promising all the way up to the very last step, at which the proof either stands or falls.

As such, it would be extremely useful to know ahead of time how difficult or easy a given congruence family will be, before one has made a serious committment.  To this end alone, this classification is beneficial.

The second advantage of course is that it gives us some sense of where the most important or promising breakthroughs in the subject will lie.  The comparatively easy congruence families are proved by methods which at this point are essentially routine.  The most difficult families on the other hand are either standing conjectures, or proved with methods which are not yet wholly systematized and which cannot be applied generally.  A good classification system can encourage the community to focus its efforts on the most difficult problems.

\section{General Method}

One noteworthy observation is that these families tend to occur for the Fourier coefficients of modular forms.  For example, the multiplicative inverse of the Dedekind eta function has the expansion

\begin{align*}
\frac{1}{\eta(\tau)} = \frac{1}{q^{1/24}(q;q)_{\infty}} = \sum_{n=0}^{\infty} p(n) q^{n-1/24}.
\end{align*}  Similarly, the counting function for partitions into distinct parts is

\begin{align*}
\frac{\eta(2\tau)}{\eta(\tau)}=q^{1/24}\frac{(q^2;q^2)_{\infty}}{(q;q)_{\infty}} = q^{1/24}(-q;q)_{\infty}= \sum_{n=0}^{\infty} p_D(n) q^{n+1/24}.
\end{align*}  It is natural to suspect that something about modularity might play a role in the occurrence of such congruence families.  What is still more difficult to explain is that this property, which superficially has the same form for so many different partition functions, varies significantly in its underlying complexity.  Our initial example for $p(n)$ already makes this clear: the case for $\ell=5,7$ is comparatively straightforward.  On the other hand, the case for $\ell=11$ is not easy at all.

A more contemporary case is the Andrews--Sellers congruence family for the generalized 2-color Frobenius partition function, $c\phi_2(n)$.  The family was identified in 1994 \cite{Sellers}, but it was not proved until 2012 \cite{Paule}.  Still other conjectured congruence families are resistant to proof to this day \cite{Beazer}.  Moreover, many proof techniques rely on heavy experimentation and guess work.  The proof of the Andrews--Sellers congruence family raises many questions, and relies on functions and algebraic structures which are not obvious at first sight.

A thorough review of the methods involved in proving congruence families will not be given here.  Useful references include \cite[Chapters 7-8]{Knopp} for the classical theory.  Some expositions of the methods behind more recent work can be found in \cite{Paule} and \cite{Baner5}.  The purpose here is simply to give a brief outline of the techniques involved, and the various complications that arise from the topology of the associated modular curve.

Suppose that one wishes to prove a family of the form
\begin{align}
a(n)\equiv 0\pmod{\ell^{\alpha}}&\text{ whenever }\Lambda n\equiv 1\pmod{\ell^{\alpha}},\label{cngfamdefn}
\end{align} in which the terms $a(n)$ are the integer coefficients of a function
\begin{align}
\mathcal{F}(\tau) := \prod_{\delta| M} (q^{\delta};q^{\delta})_{\infty}^{r_{\delta}},\label{coefdefn}
\end{align} with $q=e^{2\pi i\tau}$ and $\tau\in\mathbb{H}$, $M$ is some positive integer, and $(r_{\delta})_{\delta| M}$ is an integer-valued vector.  This is of course, up to an exponential factor, a modular form---indeed, an eta quotient.

The standard approach is to construct a sequence $\mathcal{L}:=\left( L_{\alpha} \right)_{\alpha\ge 1}$ in which
\begin{align}
L_{\alpha} = \phi(\alpha)\sum_{\Lambda n\equiv 1\bmod{\ell^{\alpha}}} a(n)q^{\left\lfloor n/\ell^{\alpha}\right\rfloor}.\label{Lalphadefn}
\end{align}  We choose the prefactor $\phi(\alpha)$ to be an integer power series with leading coefficient 1 (so that no prime divides it), with the property that $L_{\alpha}$ is a modular function over the congruence subgroup $\Gamma_0(N)$ for some positive integer $N$ when we take $q=e^{2\pi i\tau}$ with $\tau\in\mathbb{H}$.  Equivalently \cite[Chapter VI, Theorem 4A]{Lehner}, $L_{\alpha}$ corresponds to a meromorphic function on the classical modular curve $\mathrm{X}_0(N)$ with possible poles only at the cusps (the compactification points of the curve).

It is important to note that for a congruence family defined by (\ref{cngfamdefn})-(\ref{coefdefn}), there will exist a congruence subgroup $\Gamma_0(N)$ for which such modular functions (\ref{Lalphadefn}) exist, and we know how to compute it \cite{Radu}, \cite{Smoot1}.

Next, we construct a sequence (often an alternating pair) of linear operators $U^{(\alpha)}$ such that
\begin{align}
U^{(\alpha)}\left( L_{\alpha} \right) = L_{\alpha+1}.
\end{align}  These are generally modified Atkin--Lehner $U_{\ell}$ operators.

The idea is to demonstrate that $L_{\alpha}$ is divisible by $\ell^{\alpha}$ for every $\alpha\ge 1$.  We can take advantage of certain powerful finiteness conditions imposed by modularity, and express $L_{\alpha}$ in terms of some convenient reference functions over which we have greater control.

For example, because $L_{\alpha}$ may be considered as a meromorphic function on $\mathrm{X}_0(N)$, a compact Riemann surface, we have access to the following powerful result:

\begin{theorem}\label{riemannsurfacetheorema}
Let $\mathrm{X}$ be a compact Riemann surface, and let $\hat{f}:\mathrm{X}\longrightarrow\mathbb{C}$ be analytic on all of $\mathrm{X}$.  Then $\hat{f}$ must be a constant function.
\end{theorem}

\begin{corollary}
For a given $N\in\mathbb{Z}_{\ge 1}$, if $f\in\mathcal{M}\left(\Gamma_0(N)\right)$ has positive order at every cusp of $\Gamma_0(N)$, then $f$ must be a constant.
\end{corollary}

In other words, because $L_{\alpha}$ is not a constant, it must have a pole on at least one cusp of the curve.  On the other hand, we can multiply $L_{\alpha}$ by some eta quotient $z$ which annihilates every pole that $L_{\alpha}$ has except a single cusp.  For any rational number $r$, we will denote $[r]_N$ as the cusp of $\mathrm{X}_0(N)$ which contains modular transformations applied to $r$.  We will almost exclusively focus on $[0]_N$.  Part of the reason for this is simple convention.  One advantage is that it is relatively easy to study a function's behavior at $[0]_N$ and $[\infty]_N$.  We choose $[0]_N$ rather than $[\infty]_N$ because our function sequences will generally have decreasing order at our cusp of choice.  Selecting $[\infty]_N$ would give our functions a relatively cumbersome principal part.  If instead we focus on $[0]_N$, we can generally control or restrict the order at $[\infty]_N$ to something much better behaved.

Pushing $L_{\alpha}$ to the space of functions with a pole at a single point is useful, since we can generally describe the entire space of functions that live at a single cusp (e.g., $[0]_N$) in a straightforward manner dependent on the topology of $\mathrm{X}_0(N)$.

Henceforth, we denote the space of such functions as
\begin{align}
\mathcal{M}^{0}\left( \mathrm{X}_0(N) \right).
\end{align}  Our first objective is to understand what $\mathcal{M}^{0}\left( \mathrm{X}_0(N) \right)$ looks like.

\section{Genus}

Theorem \ref{riemannsurfacetheorema} gave us an important restriction on the behavior of meromorphic functions on $\mathrm{X}_0(N)$.  We now examine a different and equally important restriction in the form of the Lückensatz, or gap theorem, of Weierstrass:

\begin{theorem}[Weierstrass]\label{WGT}
Let $\mathrm{X}$ be a compact Riemann surface, and let 
\begin{align*}
f:\mathrm{X}\longrightarrow\hat{\mathbb{C}}
\end{align*} be holomorphic over $\mathrm{X}$, except for a pole at a point $p\in\mathrm{X}$.  Then the order of $f$ at $p$ can assume any negative integer, with exactly $\mathfrak{g}\left( \mathrm{X} \right)$ exceptions, which must be members of the set $\{1, 2, ..., 2\mathfrak{g}-1\}$.
\end{theorem}

See \cite{Paule2} for a proof.  As one sees in the previous section, we generally use as reference functions the generators of the space of functions which live at a given cusp of $\mathrm{X}_0(N)$ (i.e., those functions which are holomorphic everywhere on $\mathrm{X}_0(N)$ except a given cusp).

\subsection{Genus 0}

When the genus is 0 the problem is especially easy, since the order of a function which lives at a cusp can have any order.  This means that we can find a function $x$ which is holomorphic along all of $\mathrm{X}_0(N)$ except for $[0]_N$, at which $x$ has order $-1$.  It is clear, then, that $x^n$ has order $-n$ at $[0]_N$.

Supposing that $x(-1/\tau)^n$ has the expansion
\begin{align*}
x(-1/\tau)^n = \frac{b(-n)}{q^n} + \frac{b(-n+1)}{q^{n-1}} +\frac{b(-n+2)}{q^{n-2}} + ... + \frac{b(-1)}{q} + b(0) + ...,
\end{align*} and that we have a function $f$ which lives at $[0]_N$, having principal part and constant, say,
\begin{align*}
f(-1/\tau) = \frac{a(-n)}{q^n} + \frac{a(-n+1)}{q^{n-1}} +\frac{a(-n+2)}{q^{n-2}} + ... + \frac{a(-1)}{q} + a(0) + ...,
\end{align*} then of course we may have
\begin{align*}
f(-1/\tau) - \frac{a(-n)}{b(n)}x(-1/\tau)^n = \frac{c(-n+1)}{q^{n-1}} +\frac{c(-n+2)}{q^{n-2}} + ... + \frac{c(-1)}{q} + c(0) + ...,
\end{align*} for some complex numbers $c(k)$.  Of course, we may repeat this process for $x^{n-1}$, and so on until we have completely removed the principal part of $f(-1/\tau)$.  The resulting expression
\begin{align*}
f(-1/\tau) - \left( d(n)x(-1/\tau)^n + d(n-1)x(-1/\tau)^{n-1} + ... + d(1)x(-1/\tau) + d(0)\right),
\end{align*} is a modular function which is holomorphic along the whole of $\mathrm{X}_0(N)$ with no poles at all, and must therefore be a constant by Theorem \ref{riemannsurfacetheorema}.  Since we have also reduced the constant, we must have equality to 0.  We may of course rearrange terms and return to $\tau$ rather than $-1/\tau$:
\begin{align*}
f = d(n)x^n + d(n-1)x^{n-1} + ... + d(1)x + d(0).
\end{align*}  Therefore, any function which lives at $[0]_N$ is necessarily a polynomial in $x$, and we have
\begin{align}
\mathcal{M}^{0}\left( \mathrm{X}_0(N) \right) = \mathbb{C}[x].
\end{align}  Indeed, for most classical congruence families, including the first two cases of Theorem \ref{ramp5711}, the associated functions $L_{\alpha}$ in (\ref{Lalphadefn}) are members of $\mathcal{M}^{0}\left( \mathrm{X}_0(N) \right)$, in which the associated $\mathrm{X}_0(N)$ has genus 0.  In this case we have the straightforward problem of expressing $L_{\alpha}$ as a polynomial in $x$, and determining how to properly express $U^{(\alpha)}(x^m)$.

Even in more contemporary cases when the cusp count is large, one is fundamentally dealing with manipulation of a single function.  As such, the ``bookkeeping" element of the problem is relatively minimal.

\subsection{Genus 1}

The first nontrivial example of the gap theorem emerges when the genus is 1.  In this case, it is clear that the only impossible order for a function living at $[0]_N$ is $-1$, since if such a function $f$ existed, we could achieve any order at $[0]_N$ simply by taking powers $f^n$.  In particular, there ought to be functions $x,y$ with orders $-2, -3$ (respectively) at $[0]_N$.  This gives us
\begin{align}
\mathcal{M}^{0}\left( \mathrm{X}_0(N) \right) = \mathbb{C}[x]+y\mathbb{C}[x].
\end{align}  One immediately sees the complications of the genus.  Suppose that the $L_{\alpha}$ in (\ref{Lalphadefn}) are still within$\mathcal{M}^{0}\left( \mathrm{X}_0(N) \right)$.  Now, instead of worrying only about $U^{(\alpha)}(x^m)$ for a certain number of consecutive values of $m$, now we also have to concern ourselves with $U^{(\alpha)}(yx^m)$ for just as many values of $m$.  Thus in going from genus 0 to genus 1 we have effectively doubled the amount of work that must be done.

In contrast to the first two cases of Theorem \ref{ramp5711}, the third case is a family modulo powers of 11 and is associated with $\mathrm{X}_0(11)$, a curve of cusp count 2 and genus 1.  Watson proved Theorem \ref{ramp5711} for $\ell=5,7$, but he did not produce a proof of the case $\ell=11$ for Theorem \ref{ramp5711}, which he predicted could be characterized as ``boring" \cite{Watson}.  After producing a proof of this case, Atkin came to agree with Watson's prediction \cite{Atkin}.  While the word ``boring" is perhaps disputable, nevertheless one can observe that the proof method is not substantially different from the proofs in the cases $\ell=5,7$.  The most important complication of the problem lies in the doubling of work done by increasing the genus from 0 to 1.  This, combined with the relatively large size of 11, serves to make the proof very tedious indeed.

A more contemporary example occurred recently, when Banerjee and I proved a set of congruences for $d_2(n)$ by powers of 7 \cite{Baner7}.  The associated modular curve was $\mathrm{X}_0(14)$, a curve of cusp count 4 and genus 1.  We were surprised to discover that the techniques of localization described above are indeed sufficient to complete the proof.  On the other hand, localization fails to resolve the Andrews--Sellers congruences associated with $\mathrm{X}_0(20)$, a curve of cusp count 6 and genus 1.

It seems clear that the cusp count, not the genus, is responsible for the genuine difficulty of the associated congruence family.  The role of the genus, in comparison, is to make the family more ``tedious" to prove.  Because the rank of the assoicated $\mathbb{C}[x]$ module grows as the genus grows, we find ourselves with many more placeholding functions, and a great deal of bookkeeping must be done to properly keep track of what is happening when a given Atkin--Lehner operator is applied.

\subsection{Genus $\ge 2$}

As the genus grows, this problem of increasing complexity continues to apparently get worse.  This has been commented on by Gordon and Hughes \cite{Gordon2}.  

Regarding the question of general representation of the functions of $\mathcal{M}^{0}\left( \mathrm{X}_0(N) \right)$, it must be remarked that there are other intriguing possibilities, e.g. \cite[Theorem 2.1.5]{Gannon}; nevertheless, for now we will stick with a straightforward extrapolation of the cases of genus 0, 1 that we have already examined.  We will suppose that

\begin{align}
\mathcal{M}^{0}\left( \mathrm{X}_0(N) \right) = \sum_{k=0}^{v} y_k\mathbb{C}[x],
\end{align} for some $v\ge 1$ and functions $x$ and $y_0, y_1,..., y_{v}$, with $y_0=1$.

What we have, then, is the expression of our relevant space of functions as a $(v+1)$ rank $\mathbb{C}[x]$-module, for some $v\in\mathbb{Z}$, and some key reference function $x$.  We want to determine the relationship that $L_{\alpha}$ has to this module.  Our problem, then, is largely (but not entirely) redirected to the problem of studying the effects of $U^{(\alpha)}$ on various combinations of the functions $y_k,$ $x^m$.  At least in principle, this is a complication, but not a problem of accessibility.

\section{Cusp Count}

So much for the effects of the genus on our functions of interest.  We now have to examine what is possibly a more important number: the number of compactification points on $\mathrm{X}_0(N)$.

Taking $\varphi(n)$ as Euler's totient function, we have \cite[Chapter 3, Sections 3.1, 3.8]{Diamond}
\begin{align}
\epsilon_{\infty}\left( \mathrm{X}_0(N) \right) = \sum_{\delta|N}\varphi\left( \mathrm{gcd}\left(\delta,N/\delta\right) \right).
\end{align}  Notice that $\epsilon_{\infty}\left( \mathrm{X}_0(1) \right)=1$.  For $N$ a prime number we will have
\begin{align}
\epsilon_{\infty}\left( \mathrm{X}_0(N) \right) = \varphi\left( \mathrm{gcd}\left(1,N\right) \right) + \varphi\left( \mathrm{gcd}\left(N,1\right) \right) = 2\varphi(1)=2.
\end{align}  For $N$ composite we have
\begin{align}
\epsilon_{\infty}\left( \mathrm{X}_0(N) \right) = 2\varphi\left( 1 \right) + \sum_{\substack{\delta|N,\\ 1<\delta\le\sqrt{N}}} 2\varphi\left( \mathrm{gcd}\left(\delta,N/\delta\right) \right) > 2.
\end{align}  Thus, the cusp count will be 2 if and only if $N$ is prime.  Moreover, for a composite non-square, we have

\begin{align}
\epsilon_{\infty}\left( \mathrm{X}_0(N) \right) = \sum_{\substack{\delta|N,\\ \delta<\sqrt{N}}} 2\varphi\left( \mathrm{gcd}\left(\delta,N/\delta\right) \right).
\end{align}  For $N$ a square, we have

\begin{align}
\epsilon_{\infty}\left( \mathrm{X}_0(N) \right) = \sum_{\substack{\delta|N,\\ \delta<\sqrt{N}}} 2\varphi\left( \mathrm{gcd}\left(\delta,N/\delta\right) \right) + \varphi\left(\sqrt{N}\right).
\end{align}  Finally, notice from the multiplicative properties of the totient function that $\varphi\left(\sqrt{N}\right)$ is even except when $N=1$ and $N=4$.  We therefore have the following result:

\begin{lemma}
The classical modular curve $\mathrm{X}_0(N)$ has even cusp count with two exceptions for $N=1,4$.  The cusp count is exactly 2 if and only if $N$ is prime.
\end{lemma}

\subsection{Prime $N$}

For prime $N$, the problem turns out to be a very straightforward one.  One can construct $\phi(\alpha)$ and $U^{(\alpha)}$ together so that for every $\alpha\ge 1$, $L_{\alpha}$ has positive order at the cusp $[\infty]$, i.e., $L_{\alpha}$ vanishes at $q=0$.  But we still have Theorem \ref{riemannsurfacetheorema}: $L_{\alpha}$ has to have a pole somewhere.  But when $N$ is prime, there are only two cusps: $[\infty]$ and $[0]_N$.  Since the cusps are the only places where $L_{\alpha}$ can admit a pole, we are forced to accept that
\begin{align}
L_{\alpha}\in\mathcal{M}^{0}\left( \mathrm{X}_0(N) \right).
\end{align} We are thus reduced to the problem of expressing $L_{\alpha}$ in terms of these reference functions, and determining the effect of $U^{(\alpha)}$ on monomials of the form $y_kx^m$.

The classic strategy, then, is to express
\begin{align}
\frac{1}{\ell^{\alpha}}L_{\alpha} = \sum_{1\le k\le v}\sum_{m\ge 0} s_{k,\alpha}(m)\ell^{\theta(k,\alpha,m)}y_k x^m,
\end{align} in which $s_{k,m}$ is integer-valued and has finite support in $m$.  Applying $U^{(\alpha)}$, we have
\begin{align}
\frac{1}{\ell^{\alpha}}L_{\alpha+1} = \frac{1}{\ell^{\beta}}U^{(\alpha)}\left(L_{\alpha}\right) = \sum_{1\le k\le v}\sum_{m\ge 0} s_{k,\alpha}(m)\ell^{\theta(k,\alpha,m)}U^{(\alpha)}\left(y_k x^m\right).
\end{align}  We use a modular equation to construct a recurrence in how $U^{(\alpha)}$ affects powers of $x$.  We usually find that there exists an integer-valued function $h_{j,\alpha}(m,r)$ which has finite support in $r$ such that
\begin{align}
U^{(\alpha)}\left(y_k x^m\right) = \sum_{1\le j\le v}\sum_{r\ge 0} h_{j,\alpha}(m,r)\ell^{\pi(j,\alpha,m,r)}y_j x^r.
\end{align}  Combining the latter two expressions, we have
\begin{align}
\frac{1}{\ell^{\alpha}}L_{\alpha+1} = \sum_{\substack{1\le k\le v,\\ 1\le j\le v,\\ m\ge 0,\\ r\ge 0}}s_{k,\alpha}(m)h_{j,\alpha}(m,r)\ell^{\theta(k,\alpha,m)+\pi(j,\alpha,m,r)}y_j x^r.
\end{align}  One then tries to show that
\begin{align}
\theta(k,\alpha,m)+\pi(j,\alpha,m,r)\ge\theta(k,\alpha+1,r)+1.
\end{align}  This allows us to complete an induction argument.

This is the most classical approach; it is very well-understood, even by Watson \cite{Watson} and Ramanujan \cite{BO}.

\subsection{Composite $N$: Cusp Count 4}

The strategy for prime $N$ is especially simple in principle; indeed, when the genus is also 0, the Riemann surface structure is so simple that it need not be directly consulted at all.

Notice that we immediately lose this advantage when $N$ becomes composite.  At this point the cusp count function will give us more than two cusps.  As such, even if we can adjust $\phi(\alpha)$ so that each $L_{\alpha}$ has positive order at $[\infty]$, we are left with two cusps in the case that $N=4$, and at least three cusps for any other composite number.

However, we can resolve this in the following clever way: we can construct an eta quotient $z$ such that
\begin{align}
z^nL_{\alpha}\in\mathcal{M}^{0}\left( \mathrm{X}_0(N) \right)=\sum_{k=0}^{v} y_k\mathbb{C}[x],
\end{align} with $n$ some positive integer depending on $\alpha$.  Notice that we can in fact construct $z$ such that it has positive order at \textit{every} cusp besides $[0]_N$ \cite[Lemma~20]{Radu}, \cite{Newman}, so that $z$ can indeed be chosen to push $L_{\alpha}$ into the space of functions at $[0]_N$ for every $\alpha\ge 1$.

We of course arrange so that $z\in\mathcal{M}^{0}\left( \mathrm{X}_0(N) \right)$.  Indeed, we often preferentially choose a $z$ that can be constructed strictly in terms of our function $x$.  If we define
\begin{align}
\mathcal{S} :=\left\{ z^n : n\ge 0 \right\},
\end{align} then we necessarily have
\begin{align}
L_{\alpha}\in\sum_{k=0}^{v} y_k\mathbb{C}[x]_{\mathcal{S}};
\end{align} thus we can express $L_{\alpha}$ as a rational polynomial with restrictions on the denominator.  It is the manipulation of such rational function expressions of $L_{\alpha}$ that embodies the method of proof by localization.

The strategy, then, is a natural extension of the case in which $N=2$.  Suppose that
\begin{align}
\frac{1}{\ell^{\alpha}}L_{\alpha} = \frac{1}{z^n}\sum_{1\le k\le v}\sum_{m\ge 0} s_{k,\alpha}(m)\ell^{\theta(k,\alpha,m)}y_k x^m,
\end{align} in which $s_{k,m}$ is integer-valued and has finite support in $m$.  Applying $U^{(\alpha)}$, one determines that 
\begin{align}
\frac{1}{\ell^{\alpha}}L_{\alpha+1} = \frac{1}{\ell^{\beta}}U^{(\alpha)}\left(L_{\alpha}\right) = \sum_{1\le k\le v}\sum_{m\ge 0} s_{k,\alpha}(m)\ell^{\theta(k,\alpha,m)}U^{(\alpha)}\left(y_k x^m/z^n\right).
\end{align}  As in the case for $N$ a prime, we use a modular equation recurrence in $x$ and $z$ to express
\begin{align}
U^{(\alpha)}\left(y_k x^m/z^n\right) = \frac{1}{z^{n'}}\sum_{1\le j\le v}\sum_{r\ge 0} h_{j,\alpha}(m,n,r)\ell^{\pi(j,\alpha,m,n,r)}y_j x^r,
\end{align} with $n'$ some positive integer dependent on $k,m,n$, and $h_{j,\alpha}(m,n,r)$ an integer-valued function with finite support in $r$.  Combining these two results, we have
\begin{align}
\frac{1}{\ell^{\alpha}}L_{\alpha+1} = \frac{1}{z^{n'}}\sum_{\substack{1\le k\le v,\\ 1\le j\le v,\\ m\ge 0,\\ r\ge 0}}s_{k,\alpha}(m)h_{j,\alpha}(m,n,r)\ell^{\theta(k,\alpha,m)+\pi(j,\alpha,m,n,r)}y_j x^r.
\end{align}  We finally try to show that
\begin{align}
\theta(k,\alpha,m)+\pi(j,\alpha,m,n,r)\ge\theta(k,\alpha+1,r)+1.\label{localpower}
\end{align}  One important complication is that unlike the case for prime $N$, the condition (\ref{localpower}) will sometimes fail.  One must therefore study the functions $h_{j,\alpha}$ to properly compensate for this failure.  Examples of overcoming this difficulty can be found in \cite{Smoot}, \cite{Smoot0}, and in a more careful manner involving the concept of the congruence kernel in \cite{Baner5}, or in the congruence ideal sequence in \cite{Baner7}.

\subsection{Composite $N$: Cusp Count 6}

Notice that the strategy we gave in the case of cusp count 4 is a natural extension of our methods for cusp count 2 (i.e., for prime $N$).  All indications are that these methods are completely general for cusp count 2 or 4.

These techniques described above do not at present work with congruence families associated with cusp count 6.  Perhaps the most important example of such a family is the Andrews--Sellers congruence family \cite{Paule}.  Others include a set of congruences for the Rogers--Ramanujan subpartition function \cite{Smoot2}, a new family discovered by Sellers and proved recently \cite{SellersA}, and a set of standing conjectures by Beazer for singular moduli for powers of various primes \cite{Beazer}.

I have made multiple attempts to prove the Andrews--Sellers congruences via localization.  These attempts have failed thus far.  Indeed, the techniques used by Paule and Radu in \cite{Paule} are altogether different, relying on Hauptmoduln over modular curves of prime level dividing the proper level.  Similar techniques have been employed, e.g., in \cite{Smoot2}, and other more recent methods have been explored in \cite{SellersA}.

\subsection{Composite $N$: Odd Cusp Count}

As we have seen, $\epsilon_{\infty}\left( \mathrm{X}_0(N) \right)$ is odd only for $N=1,4$.  Now, $N=1$ is perhaps trivial, at least for our purposes.  However, $N=4$ is a more interesting case.  In general, for a congruence family mod powers of a prime $\ell$, we must have $\ell|N$ with $N$ the level of the associated modular curve.  As such, any congruence family for a set of Fourier coefficients associated with $\mathrm{X}_0(4)$ ought to be modulo powers of 2.

At present we are unaware of any such congruence family, although it would be very interesting to know if one exists.

\section{Classification}

We claim that there are two significant qualities to consider with a congruence family.  The first is ``difficulty"---this is a measure of the accessibility of the family to a proof.  The second is ``tediousness"---this is a measure of the bookkeeping necessary to properly produce a proof.  With these two numbers, one can properly classify congruence families, at least when they occur for the Fourier coefficients of eta quotients.

The two topological numbers which are important to us are the cusp count and the genus.  Both of these numbers are important, but by and large, $\epsilon_{\infty}$ is the more important number than $\mathfrak{g}$.

\subsection{Summary}

We organize our classification in the following table.  The bottom row gives the genus of the associated modular curve, while the leftmost column gives the cusp count.  When both numbers are small, the associated congruence family may be understood through classical techniques.  As either number grows, the difficulty in proving the associated family grows.  By and large, the cusp count appears the greater difficulty, due to the difficulties associated with the congruence kernel.  However, the genus can create substantial algebraic difficulties on its own, and neither class of difficulties should be underestimated.

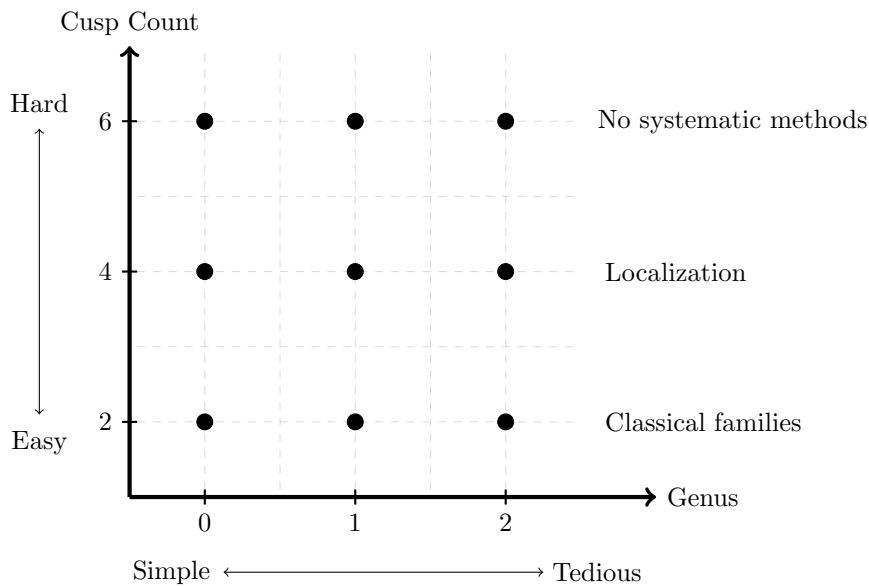
\begin{center}
\begin{tikzpicture}
\draw[help lines, color=gray!30, dashed] (-2.9,-1.9) grid (2.9,3.9);
\draw[->,ultra thick] (-3,-2)--(4,-2) node[right]{Genus};
\draw[->,ultra thick] (-3,-2)--(-3,4) node[above]{Cusp Count};
\filldraw[black] (-2,-1) circle (3pt);
\filldraw[black] (0,-1) circle (3pt);
\filldraw[black] (2,-1) circle (3pt);
\filldraw[black] (-2,1) circle (3pt);
\filldraw[black] (0,1) circle (3pt);
\filldraw[black] (2,1) circle (3pt);
\filldraw[black] (-2,3) circle (3pt);
\filldraw[black] (0,3) circle (3pt);
\filldraw[black] (2,3) circle (3pt);
\draw[-][thick] (-2.9,-1) -- (-3.1,-1) node[left]{2};
\draw[-][thick] (-2.9,1) -- (-3.1,1) node[left]{4};
\draw[-][thick] (-2.9,3) -- (-3.1,3) node[left]{6};
\draw[-][thick] (-2,-1.9) -- (-2,-2.1) node[below]{0};
\draw[-][thick] (0,-1.9) -- (0,-2.1) node[below]{1};
\draw[-][thick] (2,-1.9) -- (2,-2.1) node[below]{2};
\draw (3.2,-1) node[right]{Classical families};
\draw (3.2,1) node[right]{Localization};
\draw (3.1,3) node[right]{No systematic methods};
\draw[<->] (-1.75,-3) -- (2.5,-3) node[right]{Tedious};
\draw (-1.8,-3) node[left]{Simple};
\draw (-4.2,-1) node[below]{Easy};
\draw (-4.2,3) node[above]{Hard};
\draw[<->] (-4.2,-0.9) -- (-4.2,2.9);
\end{tikzpicture}
\end{center}

The vertical axis represents a measure of ``difficulty" in proving a given congruence family.  The horizontal axis represents a measure of ``tediousness."

Ramanujan's classic families for $p(n)$ modulo powers of 5, 7 arise at genus 0 and cusp count 2.  But Ramanujan's family for $p(n)$ modulo powers of 11 arises at genus 1 and cusp count 2.  This is further in the direction of ``tediousness," and gives some context to the commentary of Watson \cite{Watson} and Atkin \cite{Atkin}.

At cusp count 2, classical methods appear to suffice.  For cusp count 4, localization appears to be a suitable method.  For cusp count 6, localization in its current form does not work, and indeed we have no systematic methods in these cases.

The families which lie at cusp count 6 include the Andrews--Sellers congruence family \cite{Paule} and the open conjectures due to Beazer \cite{Beazer}.  We speculate that some congruence families of interest also lie at genus 2.  In the case that they do, we predict that the cusp count will determine overall accessibility of such families.

\subsection{Sporadic Families}

A key class of exceptions to this table takes place in the case that $N=4,8$.  The case for $N=4$ gives us the only classical modular curve with an odd cusp count (aside from the trivial case that $N=1$).  We are unaware of any congruence families associated with $N=4$, and it would be interesting to know if any exist, and whether the odd cusp count has any impact on the difficulty of the proof.

One case for $N=8$ was given in a paper cowritten with Sellers \cite{SellersB}, in which localization is unnecessary.  Another curious case is Boylan's proof \cite{Boylan} of a congruence family for traces of singular moduli modulo powers of 2 which is closely related to Beazer's conjectures in \cite{Beazer}.

In all of these cases, one is left with the peculiar sense that families corresponding to powers of 2 are somehow different from those corresponding to odd primes, and that they certainly lay outside of this classification scheme.

\section{Conclusion}

It must be emphasized that this classification system is basically one massive conjecture.  While we understand the effect that the genus has on the problem via the gap theorem, there does not appear to be any clear understanding as to \textit{why} $\epsilon_{\infty}$ has such a strong effect.  Certainly, a cusp count exceeding 2 demands a rational function approach that embodies localization; however, a cusp count 6 or greater entails greater difficulties which current systematic methods simply cannot overcome.  It is unclear why the cusps have this property.

Other questions that come to mind are as follows:

\begin{itemize}
\item Does our understanding extend to congruence families associated with genus 2 or greater?  In particular, do our systematic techniques work for a large genus when the cusp count is 4 or less?
\item Does it make more sense to consider different ways of representing the space of functions which live at $[0]_N$, for example using the approach given in \cite[Theorem 2.1.5]{Gannon}, especially for high genus?
\item Do the elliptic points of the associated curve also have some important impact on proving a given congruence family?
\item How do we explain the sporadic families that do not precisely lie on our classification table?
\item Can we determine systematic methods of proving congruence families associated with modular curves of cusp count 6 or greater?
\end{itemize}

In any case, we emphasize the utility of this classification against the potential for a ``cottage industry" to arise in the subject of partition congruences.  Certainly it appears that when the topology of a modular curve associated with a congruence family is simple, the family itself will be easy to prove.  The far more interesting and less accessible cases occur in matters of a nontrivial topology, and it is there that we recommend further work.

\section{Acknowledgments}

This research was funded in whole by the following research grants: the Austrian Science Fund (FWF): Einzelprojekte P 33933, ``Partition Congruences by the Localization Method," and the Austrian Science Fund (FWF): PAT 6428623 ``Towards a Unified Theory of Partition Congruences," \url{10.55776/PAT6428623}.  My sincerest and humblest thanks to the Austrian Government and People for their generous support.


\begin{thebibliography}{X}

\bibitem{AndrewsPaule2} G.E. Andrews, P. Paule, ``MacMahon's Partition Analysis VIII: Plane Partition Diamonds," \textit{Advances in Applied Mathematics} 27, pp. 231-242 (2001).

\bibitem{AndrewsPaule} G.E. Andrews, P. Paule, ``MacMahon's Partition Analysis XIII: Schmidt Type Partitions and Modular Forms," \textit{Journal of Number Theory} 234, pp. 95-119 (2021).

\bibitem{Atkin} A.O.L. Atkin, ``Proof of a Conjecture of Ramanujan,'' \textit{Glasgow Mathematical Journal} 8, pp. 14-32 (1967).

\bibitem{Atkin0} A.O.L. Atkin, ``Ramanujan Congruences for $p_{-k}(n)$," \textit{Canadian Journal of Mathematics} 20, pp. 168-178 (1968).

\bibitem{AtkinL} A.O.L. Atkin, J. Lehner, ``Hecke Operators on $\Gamma_0(M)$,'' \textit{Mathematische Annalen} 185, pp. 134-160 (1970).

\bibitem{Baner5} K. Banerjee, N.A. Smoot, ``The Localization Method Applied to $k$-Elongated Plane Partitions and Divisibility by 5" (submitted) (2022), \url{https://arxiv.org/abs/2208.07065}.

\bibitem{Baner7} K. Banerjee, N.A. Smoot, ``2-Elongated Plane Partitions and Powers of 7: The Localization Method Applied to a Genus 1 Congruence Family" (submitted), \url{https://arxiv.org/abs/2306.15594}

\bibitem{Beazer} M. Beazer, \textit{3-adic Properties of Hecke Traces of Singular Moduli}, Master Thesis, Brigham Young University (2021).

\bibitem{BO} B.C. Berndt, K. Ono, ``Ramanujan's Unpublished Manuscript on the Partition and Tau Functions,'' \textit{The Andrews Festschrift}, Springer-Verlag, pp. 39-110 (2001).

\bibitem{Boylan} M. Boylan, ``2-adic Properties of Hecke Traces of Singular Moduli," \textit{Mathematics Research Letters} 12 (4), pp. 593-609 (2005).

\bibitem{Diamond} F. Diamond, J. Shurman, \textit{A First Course in Modular Forms}, 4th Printing., Springer Publishing (2016).

\bibitem{Euler} L. Euler, \textit{Introductio in Analysin Infinitorum}, Chapter 16.  Marcum--Michaelem Bousquet, Lausannae (1748).

\bibitem{Gannon} T. Gannon, \textit{Moonshine Beyond the Monster: The Bridge Connecting Algebra, Modular Forms, and Physics}, Cambridge University Press, New York (2006).

\bibitem{Gordon2} B. Gordon, K. Hughes, ``Ramanujan Congruences for $q(n)$, in: M. Knopp, \textit{Analytic Number Theory.  Lecture Notes in Mathematics,} 899, Springer, Berlin, Heidelberg (1981).

\bibitem{Knopp} M. Knopp, \textit{Modular Functions in Analytic Number Theory}, 2nd Ed., AMS Chelsea Publishing (1993).

\bibitem{Lehner} J. Lehner, \textit{Discontinuous Groups and Automorphic Functions}, Mathematical Surveys and Monographs Number 8, American Mathematical Society (1964).

\bibitem{Lehner6} J. Lehner, ``Ramanujan Identities Involving the Partition Function for the Moduli $11^{\alpha}$,'' \textit{American Journal of Mathematics} 65, pp. 492-520 (1943).

\bibitem{Lehner7} J. Lehner, ``Proof of Ramanujan's Partition Congruence for the Modulus $11^3$,'' \textit{Proceedings of the American Mathematical Society} 1, pp. 172-181 (1950).

\bibitem{MacMahon} P.A. MacMahon, \textit{Combinatory Analysis}, Cambridge University Press, Cambridge, 1915-1916 (Reprinted: Chelsea, New York, 1960).

\bibitem{Newman}  M. Newman, ``Construction and Application of a Class of Modular Functions (II)," \textit{Proc. London Math. Soc.}, 3 (1959).

\bibitem{Paule}  P. Paule, S. Radu, ``The Andrews--Sellers Family of Partition Congruences," \textit{Advances in Mathematics} 230, pp. 819-838 (2012).

\bibitem{Paule2}  P. Paule, S. Radu, ``A Proof of the Weierstraß Gap Theorem not Using the Riemann--Roch Formula," \textit{Annals of Combinatorics} 23, pp. 963-1007 (2019).

\bibitem{Radu} S. Radu, ``An Algorithmic Approach to Ramanujan--Kolberg Identities," \textit{Journal of Symbolic Computation}, 68, pp. 225-253 (2015).

\bibitem{Ramanujan} S. Ramanujan, ``Some Properties of $p(n)$, the Number of Partitions of $n$", \textit{Proceedings of the Cambridge Philosophical Society} 19, pp. 207-210 (1919).

\bibitem{Rodseth} \O . R\o dseth, ``Congruence Properties of the Partition Functions $q(n)$ and $q_0(n)$,'' \textit{Arbok Univ. Bergen Mat.--Natur. Ser.} 13 (1970).

\bibitem{Sellers} J. Sellers, ``Congruences Involving F-Partition Functions,'' \textit{International Journal of Mathematics and Mathematical Sciences} 17, ppl 187-188 (1994).

\bibitem{SellersB} J.A. Sellers, N.A. Smoot, ``On the Divisibility of 7-Elongated Plane Partition Diamonds by Powers of 8," \textit{International Journal of Number Theory} Volume 20 (1), (2024).

\bibitem{SellersA} J.A. Sellers, N.A. Smoot, ``Old Meets New: Connecting Two Infinite Families of Congruences Modulo Powers of 5 for Generalized Frobenius Partition Functions" (submitted) (2024), \url{https://arxiv.org/abs/2402.18509}

\bibitem{Smoot1} N.A. Smoot, ``On the Computation of Identities Relating Partition Numbers in Arithmetic Progressions with Eta Quotients: An Implementation of Radu's Algorithm,'' \textit{Journal of Symbolic Computation} 104, pp. 276-311 (2021).

\bibitem{Smoot} N.A. Smoot, ``A Congruence Family For 2-Elongated Plane Partitions: An Application of the Localization Method,'' \textit{Journal of Number Theory} 242, pp. 112-153 (2023).

\bibitem{Smoot2} ``A Family of Congruences for Rogers--Ramanujan Subpartitions," \textit{Journal of Number Theory,} Volume 196, pp. 35-60 (2019).

\bibitem{Smoot0} N.A. Smoot, ``A Single-Variable Proof of the Omega SPT Congruence Family Over Powers of 5,'' \textit{Ramanujan Journal} 62, pp. 1-45 (2023).

\bibitem{Watson} G.N. Watson, ``Ramanujans Vermutung über Zerfallungsanzahlen," \textit{J. Reine Angew. Math.} 179, pp. 97-128 (1938).
\end{thebibliography}
\end{document}